\documentclass[11pt]{article}
\usepackage{latexsym,amsthm,verbatim,ifthen,graphicx,color,mathrsfs}
\usepackage[all]{xy}
\usepackage{color}

\usepackage{color}

\usepackage{lscape}
 \usepackage[ansinew]{inputenc}
 \usepackage{multicol}
 \usepackage[all]{xy}\xyoption{poly}
 \usepackage{mathrsfs}
\usepackage[psamsfonts]{eucal}
\usepackage{amssymb, amsmath}
\usepackage{geometry}
\usepackage{enumerate}
\usepackage{float}

\newtheorem{innercustomthm}{Theorem}
\newenvironment{customthm}[1]
  {\renewcommand\theinnercustomthm{#1}\innercustomthm}
  {\endinnercustomthm}

\newtheorem*{theorem*}{Theorem}
\newtheorem{theorem}{Theorem}[section]

\newtheorem{corollary}[theorem]{Corollary}
 \newtheorem{lemma}[theorem]{Lemma}
 \newtheorem{proposition}[theorem]{Proposition}
 \theoremstyle{definition}
 
 \theoremstyle{remark}

\def\quil{{\mathscr L}}
\def\quilc{\widehat{\mathscr L}}
\def\quic{{\mathscr C}}
\def\quip{{\mathscr P}}
\def\quicc{\widehat{\mathscr C}}
\def\calcc{\widehat{\mathcal C}}

\def\bq{\mathbb{Q}}

\def\bn{\mathbb{N}}

\def\des{{s^{-1}}}

\def\timest{\widetilde\times}

      \newcommand{\derr}{{\rm Der}}

\newcommand{\tor}{\operatorname{\text{Tor}}}

 \newcommand{\lib }{\mathbb{L}}

  \newcommand{\Hom}{\operatorname{\text{\rm Hom}}}

\newcommand{\catdgl}{\operatorname{{\bf dgl}}}

\newcommand{\catcdgl}{\operatorname{{\bf cdgl}}}
\newcommand{\catcdgle}{\operatorname{{\bf cdgl^{\Delta}}}}

  \newcommand{\Map}{\operatorname{{\rm map}}}

        \newcommand{\ad}{\operatorname{{\rm ad}}}
 \newcommand{\calc}{{\mathcal C}}

 \newcommand{\MC}{\operatorname{{\rm MC}}}

     \newcommand{\ev}{\operatorname{{\rm ev}}}
\newcommand{\mc}{{\MC}}

\newcommand{\catcdgc}{\operatorname{{\bf cdgc}}}

  \newcommand{\otimesc}{\widehat{\otimes}}

   \newcommand{\libc}{{\widehat\lib}}

\begin{document}

\title{Homology of the completion of a Lie algebra}

\author{Yves F\'elix and Aniceto Murillo}

\maketitle

\begin{abstract}
We prove that the second homology group of the completion of an infinite dimensional free Lie algebra is uncountable.
\end{abstract}

\section*{Introduction}

The Malcev correspondence \cite{mal}, which  establishes an isomorphism between the categories of rational nilpotent groups and rational nilpotent Lie algebras, associates to any Lie algebra $L$ the rational nilpotent group $\Gamma$ consisting on $L$  with the group structure given by the Baker-Campbell-Hausdorff formula. For a general rational Lie algebra $L$, its completion $\widehat{L}=\varprojlim_n L/L^n$, where $L^1= L$ and for $n>1$, $L^n = [L, L^{n-1}]$ can also be endowed of the group structure  $\Gamma$ with the BCH product. In this general context the  usual homologies of $\widehat{L}$ and of $\Gamma$ remain mysterious although they are  important from a topological point of view. Let $G$ be a finitely generated group and denote by $BG$ its classifying space. Then, see   \cite[\S12.2]{bfmt0} or \cite[\S2]{FHTII} and, there is a complete graded Lie algebra $\widehat L$ such that the group $\Gamma$ associated to $\widehat L_0$ is  the Malcev $\bq$-completion of $G$ and the fundamental group of the Bousfield-Kan $\mathbb Q$-completion of $BG$ \cite{BK}.

In \cite[Theorem 1]{IM} S. Ivanov and R. Mikhailov prove that if $\Gamma$ is the Malcev completion of a free group on two generators then $H_2(\Gamma) $ has an uncountable dimension. Here, we prove a corresponding result for the completion $\widehat L$ of a free Lie algebra with at least two generators:

\begin{customthm}{1}\label{main}
\label{tX} Let $\widehat L$ be the completion of a free Lie algebra $L$ with  $\dim L/[L,L] >1$.  Then, $H_2(\widehat L)$ is of uncountable dimension.
\end{customthm}

\section{Chains of the completion vs completion of the chains}

By default, every considered algebraic structure has $\bq$ as ground field.
Although we are inspired by the topological setting in which Lie algebras are naturally graded, we assume, unless explicitly specified otherwise, that any consider Lie algebra  is concentrated in degree 0. In this context, recall that the homology of a Lie algebra $L$ is defined as $H(L)=\tor_{UL}(\bq,\bq)$ and coincides with the homology of the {\em chain coalgebra} $\quic(L)$, see for instance \cite[p. 301]{FHT}. This  is the differential graded coalgebra $\quic(L)=(\land sL,d)$ where $s$ denotes suspension, i.e., $sL$ is concentrated in degree $1$. The differential, which vanishes on $\quic_0(L)=\bq$ and  $\quic_1(L)=sL$, is given in $\quic_p(L)$ by,
$$d (sx_1\land \dots \land sx_p)= \sum_{1\leq i<j\leq p} (-1)^{i+j}s[x_i, x_j]\land sx_1\land \dots \widehat{sx_i}\dots \widehat{sx_i}\dots \land sx_p,$$
for any $p\ge 2$. In particular, given $x,y,z\in L$,
$$\renewcommand{\arraystretch}{1.5}
\left\{
\begin{array}{l} d(sx\land sy) = -s[x,y],\\
 d (sx\land sy\land sz) = -s[x,y]\land sz + sx\land s[y,z] - sy\land s[x,z].\end{array}
 \right.
 \renewcommand{\arraystretch}{1}$$

We define $H(L)$ to be the homology of $\quic(L)$

 Now, assume that $\lib=\lib(V)$ is a free Lie algebra. Then, $\quic(\lib)$ is naturally bigraded as follows:
 $$
 \quic(\lib)=\oplus_{n\ge p\ge 0}\quic_p(\lib)(n)
 $$
 where $\quic_0(\lib)(0)=\bq$, $\quic_0(\lib)(n)=0$ for $n>0$ and, for $p\ge 1$,
 $$
 \quic_p(\lib)(n)=\oplus_{i_1+\dots+i_p=n}\,\,s\lib^{i_1}\land\dots\land s\lib^{i_p},
 $$
 being $\lib^{k}$ the sub vector space of $\lib$ generated by brackets of length $k$ in $V$.

 Observe that
\begin{equation}\label{prime}
 d\, \quic_p(\lib)(n)\subset \quic_{p-1}(\lib)(n)
 \end{equation}
 and thus, $\quic(\lib)(n)=\oplus_{p}\quic_p(\lib)(n)$ is a subcomplex of $\quic(\lib)$ for any fixed $n$. In particular we may filter $\quic(\lib)$ by the subcomplexes
 $$
 F^0=\quic(\lib)\supset F^1\supset\dots\supset F^q\supset F^{q+1}\supset\dots,
 $$
 where $F^q=\oplus_{n\ge q}\, \quic(\lib)(n)$, and consider the {\em completion of the chains on $\lib$} with respect to this filtration,
 $$
 \quicc(\lib)=\varprojlim_q \quic(\lib)/F^q = \prod_n \, \quic(\lib)(n).
 $$
Note that $\quicc_0(\lib)=\bq$ and, for $p\ge 1$, a generic element of $\quicc_p(\lib)$ can be written as a series,
$$
\sum_{n\ge p}\Phi_n\qquad\text{with}\qquad \Phi_n\in\quic_p(\lib)(n).
$$

\begin{lemma}\label{lemma1}
$$
H_p\bigl(\quicc(\lib)\bigr)=H_p(\lib)=\begin{cases} \text{$\bq\,\,$ if $p=0$},\\ \text{$V\,\,$ if $p=1$},\\   \text{$0\,\,\,$ if $p\ge 2$}.\end{cases}
$$
\end{lemma}

\begin{proof}
The second equality is a well known fact. On the other hand, in view of (\ref{prime}), the homology of $\quic(\lib)$  coincides with that of its completion.
\end{proof}

\medskip

We consider now the completion $\libc=\libc(V)$ of $\lib$ with respect to the ideals $\{\lib^{\ge q}\}_{q\ge 1}$ conforming  the central series. That is,
$$
\libc=\varprojlim_q\lib/\lib^{\ge q}.
$$
We refer to \cite[\S3.2]{bfmt0} for a detailed compendium of this particular class of Lie algebras.  An  element of this completion can be written as a series
$$
\sum_{q\ge 1}\Psi_q\qquad\text{with}\qquad  \Psi_q\in\lib^q.
$$
Therefore, for $p\ge 1$, a generic element of $\quic_p(\libc)$ is a finite sum of terms of the form,
$$
s(\sum_{i_1\ge 1}\Psi_{i_1})\land\dots\land s(\sum_{i_p\ge 1}\Psi_{i_p}),\qquad \Psi_{i_j}\in\lib^{i_j}.
$$
Write formally this term as
$$
\sum_{n\ge p}\quad \sum_{i_1+\dots+i_p=n} s\Psi_{i_1}\land\dots\land s\Psi_{i_p}
$$
and observe that, since $s\sum_{i_1+\dots+i_p=n} \Psi_{i_1}\land\dots\land s\Psi_{i_p}\in\quic_p(\lib)(n)$,
the above term is a well defined element in $\quicc(\lib)$. This exhibits the chains of the completion of $\lib$ as a subcomplex of the completion of the chains on $\lib$,
$$
\quic(\libc)\subset \quicc(\lib).
$$
Our main result illustrates how different  these complexes are even from the homological point of view.

\section{Chains of occurrence 2}

In this section and in the next one,  $\lib=\lib(a,b)$ will denote the free Lie algebra on two generators $a$ and $b$.

In this case we  define the {\em occurrence } (of $b$) of an iterated bracket of $\lib$ to be the number of  times that $b$ appears in the bracket. An element in $\lib$ (respec. $\libc$) is of {\em occurrence $k$} if it is a sum (respec. series) of iterated brackets of occurrence $k$. This definition extends to both $\quic(\libc)$ and $\quicc(\lib)$:

By definition each scalar in $\bq$ has occurrence $0$. For $p\ge 1$, given elements $x_1, \dots,x_p $ of $\lib$ or $\libc$, of respective occurrence $k_1, \dots , k_p$, the {\em occurrence} of the element $sx_1\land \dots \land sx_p$ in $\quic_p(\lib)$ or $\quic_p(\libc)$ is $k_1+\dots + k_p$. This already defines occurrence in $\quic(\libc)$. Finally, an element $\Phi=\sum_{n\ge 1}\Phi_n\in\quicc_p(\lib)$ has occurrence $k$ it each $\Phi_n$ does.

Observe that the differential in both $\quicc(\lib)$ and $\quic(\libc)$ respects occurrence $k$ as long as $k\ge 2$. We then focus on  $k=2$ and define
$$
\calcc(\lib)\subset \quicc(\lib)\qquad\text{and}\qquad  \calc(\libc)\subset \quic(\libc)
$$
as the subcomplexes consisting on elements of  occurrence $2$. The following properties of these complexes will be used:

\begin{proposition}\label{propo1}
\begin{enumerate}
\item[(i)]  $\calcc(\lib)=\calcc_{\le 3}(\lib)$ and $\calc(\libc)=\calc_{\le 3}(\libc)$.
\item[(ii)] $H_p\bigl(\calcc(\lib)\bigr)=0$, for $p\ge 2$.
\item[(iii)]  $\calcc_3(\lib)$ is isomorphic to the vector space of  infinite antisymmetric square matrices with coefficients in $\mathbb Q$.
\item[(iv)]  $\calc_3(\libc)$ is isomorphic to the sub vector space generated by the antisymmetric matrices of the form
$CD^t-DC^t$
with $C,D$ infinite 1-column matrices.
\item[(v)] $\calc(\libc)$ is a retract of  $\quic(\libc)$ and therefore $H(\calc(\libc))$ is a retract of $H(\libc)$.
\end{enumerate}
\end{proposition}

In what follows we strongly use the following notation:
 for any integer $r\geq 0$, we denote by $b_{r}$ the bracket of length $r+1$ and occurrence $1$ in $\lib$ given by
$$b_{r} = \bigl[[\dots \bigl[[b,a],a\bigr]\dots ,a\bigr]=(-1)^r\ad^r_a(b),$$
being $\ad$ the adjoint operator.
\begin{proof} (i) Simply observe that for any $p\ge 4$, any element $sx_1\land\dots\land sx_p\in\quic_p(\lib)$ of occurrence $2$ necessarily vanishes since at least two of the $x_i$'s are multiple of $a$.

(ii) Trivial in view of Lemma \ref{lemma1} and the fact that the differential respects occurrence.

(iii) For the sake of simplicity we drop the suspension symbol in the formulae concerning chains. Observe that every element $\Phi\in\calcc_3(\lib)$ can be uniquely written as
$$
\Phi=\sum_{n\ge 3}\Phi_n\qquad\text{with} \qquad \Phi_n= \sum\limits_{\substack{p+q=n-3\\ p>q\ge 0}}\lambda_{pq}\,\,b_{p}\land b_{q}\land a,\quad\lambda_{pq}\in\bq.
$$
This is naturally identified with the infinite dimensional antisymmetric square matrix
\begin{equation}\label{matrix}
M(\Phi)=(a_{ij})_{i,j\ge 0}\qquad\text{with}\quad a_{ij}=\lambda_{ij}\quad\text{for}\quad i>j.
\end{equation}

(iv) An element of $\calc_3(\libc)$ is a linear combination of terms of the form
$$
(\sum_{p\ge 0}\mu_p\, b_{p})\land (\sum_{q\ge 0}\nu_q\, b_{q})\land a
$$
Let $C$ and $D$ the infinite 1-column  matrices given by the sequences $\{\mu_p\}_{p\ge 0}$ and $\{\nu_q\}_{q\ge 0}$ respectively. Then, an elementary computation shows that the above term can be expressed as
$$
\sum_{n\ge 3}\,\,\, \sum\limits_{\substack{p+q=n-3\\ p>q\ge 0}}\lambda_{pq}\,\,b_{p}\land b_{q}\land a,
$$
being $\lambda_{pq}$ the corresponding coefficient of the infinite dimensional antisymmetric matrix $CD^t-DC^t$.

(v) Indeed $\quic(\libc)= \calc(\libc)\oplus A$ where $A$ is the subcomplex generated by elements of occurrence $\neq 2$.
\end{proof}

\section{The proof}
 In this section we  prove Theorem \ref{main} for the completion of a free Lie algebra on two generators. This is a direct consequence of Theorem 1' below.

 \vspace{3mm}\noindent {\bf Theorem 1'.} \emph{ With the above notations, the dimensions of $H_2(\calc(\libc)$ and $H_2(\libc)$ are uncountable.}

 \vspace{3mm}  Again, the notation of the chains is simplified by omitting the suspension sign everywhere.

 For a given pair $\{r_k\}_{k\ge 1}$ and $\{s_\ell\}_{\ell\ge 1}$ of increasing sequence of integers consider the following element of $\calcc_3(\lib)$:
 $$
\alpha  =\sum_{k,\ell } \sum_{i=0}^{s_\ell-1} (-1)^i\,\,b_{r_k+i}\land b_{s_\ell-i-1}\land a.$$
 Notice that, for any integers $p,q\ge 0$ one has
 $$
 d(b_p\land b_q\land a)=b_{p+1}\land b_{q}+ b_{p}\land b_{q+1} -[b_p,b_q]\land a.
 $$
 With this in mind a straightforward computation, which includes a cancellation procedure, yields:
 $$
 d\alpha= \bigl(\sum_k b_{r_k}\bigr) \land \bigl(\sum_\ell b_{s_l}\bigr)- \omega \land a + \Bigl( \sum_{k,\ell} (-1)^{s_l-1} b_{r_k+s_l}\Bigr)\land b,$$
 where
$$\omega = \sum_{k,\ell} \,\,\sum_{i=0}^{s_\ell -1}(-1)^i [b_{r_k+i}, b_{s_\ell-i-1}].$$
Call
$$
\Omega=d\alpha$$
and observe that this is a well defined element in $\calc_2(\libc)$ which is obviously a cycle. Moreover, for any sub sequences $\{r_k\}_{k\in A}$ and $\{s_\ell\}_{\ell\in A}$, with $A\subset\bn$, the same procedure produces the cycle of $\calc_2(\libc)$
$$\Omega_A  = \bigl(\sum_{k\in A} b_{r_k}\bigr) \land \bigl(\sum_{\ell\in A} b_{s_l}\bigr) - \omega_A \land a + \left( \sum_{k,\ell\in A} (-1)^{s_l-1} b_{r_k+s_l}\right)\land b,$$
with
$$\omega_A = \sum_{k,\ell\in A} \,\,\sum_{i=0}^{s_\ell -1}\,(-1)^i \,[b_{r_k+i}, b_{s_\ell-i-1}],$$
which is the boundary of the element in $\calcc_3(\lib)$ given by
$$\alpha_A  =\sum_{k,\ell\in A }\,\, \sum_{i=0}^{s_\ell-1} \,(-1)^i \,\,b_{r_k+i}\land b_{s_\ell-i-1}\land a.$$

We plan to construct, for a particular pair of increasing sequences $(r_k)$ and $(s_\ell)$,  an uncountable number of subsequences for which the corresponding cycles $\Omega_A$ represents different homology classes in $H_2\bigl(\calc(\libc)\bigr)=H_2(\libc)$.

Let $\{r_k\}_{k\ge 1}$ and $\{s_\ell\}_{\ell\ge1}$ be sequences of integers satisfying:
$$2=s_1 <r_1 <s_2 <r_2 \dots <r_{n-1}<s_n<r_n<\dots\quad\text{with}\quad s_n>3r_{n-1} \quad\text{and}\quad r_n > 2s_n.
$$
Define
$$F = \{(r_k+i, s_{\ell}-i-1),\,\,\,  \,  0\le  i\leq s_\ell -1\,\}.$$
That is, $F$ consists of pairs $(m,n)$ for which the term $b_m\land b_n\land a$ appears as a component, possibly affected of a sign, in the element $\alpha$.

\begin{lemma}\label{lemma2}
\begin{enumerate}
\item[(i)] The elements of $F$ are  different form each other.
\item[(ii)] For $i>k$, $(r_k+s_k-s_i, s_t-1)\in F$ if and only if $t= i$.
\item[(iii)] If $(m,n)\in F$, then $(n,m)\not\in F$.
\end{enumerate}
\end{lemma}
\begin{proof}
We first show that,
\begin{equation}\label{equ0}
 k\neq p\quad\text{implies}\quad r_k+s_\ell \neq r_p+s_q\quad\text{for all}\quad \ell, q.
\end{equation}
By contradiction, suppose $r_k + s_\ell = r_p+s_q$ and $k>p$. Then  $q>\ell$. Hence, since$r_k>6r_{k-1}$ and $s_q>6s_{q-1}$, it follows that $r_k-r_p>\frac{5}{6} r_k$ and $s_q-s_\ell > \frac{5}{6} s_q$.  Therefore, $\frac{5}{6}s_q<s_q-s_\ell = r_k-r_p<r_k$, so that $q\leq k$. Then,
$$s_q-s\ell <s-q <\frac{r_k}{2} <r_k-r_p,$$
which is a contradiction.

(i) If $(r_k+i, s_\ell -i-1)=(r_{k'}+j, s_{\ell'}-j-1)$ then $r_k+s_\ell =r_{k'}+s_{\ell'}$. By (\ref{equ0}), $k=k'$ and $\ell = \ell'$

(ii) If $(r_k+s_k-s_i, s_t-1)\in F$, there exist $\ell$ and $j$ with
$$r_k+s_k-s_i = r_k+j\qquad\text{and}\qquad s_t-1 = s_\ell -j-1.$$
In this case $s_\ell -s_t = j<r_k$, and so $\ell \leq k$. When $\ell = k$, we have
$$s_k-s_i = j= s_k-s_t,$$
and so $t= i$. When $\ell<k$, the equation $s_\ell-s_t =s_k-s_i$ has no solution because $s_k >3r_{k-1}$.

(iii) Write $(a,b)= (r_t+j, s_u-j-1)$ with $j\leq s_u-1$. Suppose $(b,a)= (r_k+i, s_\ell -i-1)$. Then $s_u-r_k= i+j+1=s_\ell -r_t$, so that $s\ell +r_k= s_u+r_t$, which contradicts (\ref{equ0}).

\end{proof}

Now, for the  sub sequences let $\{r_k\}_{k\in A}$ and $\{s_\ell\}_{\ell\in A}$, with $A\subset\bn$, we denote by  $M(\alpha_A)$ the matrix associated to $\alpha_A$ by Proposition \ref{propo1} (iii). Denote then by  $M(\alpha_A)^{[k]}$ the   $(k-1)\times (k-1)$   sub matrix of $M(\alpha)$ formed by the rows
$$r_k+s_k-s_{k-1}, r_k+s_k-s_{k-2}, \dots , r_k+s_k-s_2, r_k+s_k-s_1,$$
and the columns
$$s_1-1, s_2-1, \dots , s_{k-1}-1.$$
Then, in view of (i), (ii) and (iii) of Lemma \ref{lemma2}, a simple inspection shows:

\begin{lemma}\label{lemma3}
\begin{enumerate}
\item[(i)] If $k\in A$, then $M(\alpha_A)^{[k]}$ is the anti-diagonal matrix
$$
{\raisebox{.6\height}{$\Biggl($}}{\pm 1^{{}^{\displaystyle  .^{\displaystyle .^{\displaystyle .^{\displaystyle .^{\raisebox{-.5\height}{$\displaystyle
                        \pm 1$}}}}}}}}{\raisebox{.6\height}{$\Biggr).$}}
                        $$

\item[(ii)] If $k\not\in A$, then $M(\alpha_A)^{[k]}$ is the matrix $0$.\hfill$\square$
\end{enumerate}
\end{lemma}

Next, recall that the set $\quip(\bn)$ of subsets of $\bn$, equipped with the symmetric difference operator $\Delta$, is a $\mathbb Z/2\mathbb Z$-vector space for which the subset $\quip_f(\bn)$ of finite subsets is a sub vector space. We denote by $\{A_\gamma\}_{\gamma\in \Gamma}$ a family of elements of $\quip(\bn)$ which project onto a basis of the quotient $\quip(\bn)/\quip_f(\bn)$.  Since $\quip(\bn)$ is not countable and $\quip_f(\bn)$ is a countable set, the index set $\Gamma$ is uncountable.

It follows that, given $\gamma_1, \dots , \gamma_n\in \Gamma$, the symmetric difference
$A_{\gamma_1}\Delta\dots\Delta A_{\gamma_n}
$
is an infinite set. In particular,  for some $j\in\{1,\dots,n\}$,
\begin{equation}\label{equ2}
A_{\gamma_j} \backslash (\cup_{i\neq j} A_{\gamma_i} \cap A_{\gamma_j} )
\end{equation}
is also an infinite set.

We then consider the uncountable family of cycles
\begin{equation}\label{omegas}
\Omega_{A_\gamma}\in\calc_2(\libc)\qquad\text{with}\qquad \gamma\in \Gamma,
\end{equation}
and show that they represent linearly independent classes in $H_2\bigl(\calc(\libc)\bigr)=H_2(\libc)$.
For it, we simplify the notation and write
$$
\Omega_{A_\gamma}=\Omega_\gamma\qquad\text{and}\qquad \alpha_{A_\gamma}=\alpha_\gamma.
$$

\begin{lemma}\label{lemma4} A (finite) linear combination $\sum_{\gamma}\lambda_\gamma\Omega_\gamma$ is a boundary in $\calc(\libc)$ if and only if $\sum_{\gamma}\lambda_\gamma\alpha_\gamma\in\calc_3(\libc)$.
\end{lemma}

\begin{proof} We know that $\sum_{\gamma}\lambda_\gamma\Omega_\gamma=d(\sum_{\gamma}\lambda_\gamma\alpha_\gamma)$. If we suppose  that $\sum_{\gamma}\lambda_\gamma\Omega_\gamma=d\beta$ with $\beta\in \calc_3(\libc)$ it follows that $d(\sum_{\gamma}\lambda_\gamma\alpha_\gamma-\beta)=0$. However, see (i) and (ii) of Proposition \ref{propo1}, since $\calcc_4(\lib)=0$ and $H_3\bigl(\calcc(\lib)\bigr)=0$  the differential $d$ is injective on $\calcc_3(\lib)$ and thus $\sum_{\gamma}\lambda_\gamma\alpha_\gamma=\beta\in\calc_3(\libc)$.

\end{proof}

In view of this result the proof of Theorem \ref{main} is complete when we show that, if a  linear combination of the family $\{\alpha_\gamma\}_{\gamma\in\Gamma}$  lives in $\calcc_3(\lib)$, then it must vanish.

Assume that  $\Phi= \sum_{i=1}^n \lambda_{i} \alpha_{\gamma_i}$ belongs to $\calc_3(\libc)$ and let $n$ be the minimum integer for which this happens.

On the one hand, in view of Proposition \ref{propo1}(iv), $M(\Phi)$ has finite rank since every matrix of the form $CD^t-D^tC$ is of rank at most $1$.

Now, if $n= 1$ then $\Phi=\lambda\alpha_\gamma$ for some $\gamma\in\Gamma$. By Lemma \ref{lemma3}(i) the determinant of $M(\Phi)^{[k]}=\pm\lambda$ for every $k\in A_\gamma$ and thus, the rank of $M(\Phi)$ is infinite unless $\lambda=0$.

If $n>1$, in view of (\ref{equ2}), choose $j\in \{1, \dots , n\}$ so that  $A_{\gamma_j} \backslash (\cup_{i\neq j} A_{\gamma_i} \cap A_{\gamma_j} )$ is an infinite set. Then, again by Lemma \ref{lemma3}, it follows that, for each $k\in A_{\gamma_j} \backslash (\cup_{i\neq j} A_{\gamma_i} \cap A_{\gamma_j} )$, the matrix
$M(\alpha_{\gamma_i})^{[k]}$ is invertible if $i=j$ and $0$ otherwise. Hence,  $M(\Phi)^{[k]}= M(\sum_{i=1}^n \lambda_i  \alpha_{\gamma_i})^{[k]}= \lambda_{j}  M(\alpha_{\gamma_{j}})^{[k]}$. Hence, this implies again that the rank of  $M(\Phi)$ is infinite unless $\lambda_j= 0$. By the minimality condition on $n$ we conclude that $\Phi = 0$.

\section{Consequences and related questions}

As in the past section $\lib=\lib(a,b)$ and $\libc$ denotes its completion. Easy consequences of Theorem 1'   are:

\begin{corollary}\label{cor1} Let $L'$ and $L''$ be Lie algebras and let $L$ be  the completion of the free product $L'\amalg L''$. Then, $H_2(L)$ is of uncountable dimension.
\end{corollary}

\begin{proof}  Let $a$ and $b$ be  indecomposable elements in $L'$ and $L''$ respectively. Then $\lib(a)$ and $\lib(b)$ are retract of $L'$ and  $L''$ respectively. It follows that $\libc = \mathbb L(a)\, \widehat{\amalg}\, \mathbb L(b)$ is a retract of $L$, and so $H_2(\libc)$ is a retract of $H_2(L)$. By Theorem \ref{main} the corollary follows.
\end{proof}

\vspace{3mm} In particular this implies Theorem 1 in the general case.

\begin{corollary}\label{cor2}  Let $L$ be the completion of the Lie algebra $\lib/I$ where $I$ is contained in the Lie algebra generated by Lie brackets of occurrence greater than $2$ in $b$. Then $H_2(L)$ is of uncountable dimension.
\end{corollary}

\begin{proof} Let $J$ the ideal generated by all brackets of occurrence greater than $2$. Then $I\subset J$ and, by completing the corresponding projections, we have morphisms
$$\libc\to \widehat{\lib/I} \to \widehat{\lib/J}.$$

Now denote by $\calc(\widehat{\lib/J})\subset \quic(\widehat{\lib/J})$ the injection of chains of occurrence 2, that is once again a retraction. We remark that by construction the projection
$$\calc(\libc) \to \calc(\widehat{\lib/J})$$
is an isomorphism. Therefore by Theorem 1', $H_2(L)$ is also of uncountable dimension.
\end{proof}

We now observe that analogous infiniteness properties remain even  when we restrict to   brackets of occurrence $1$.

\begin{theorem} Let $K \subset \lib$ be the ideal generated by Lie brackets of occurrence 2 in $b$.
 Then $H_2(\lib/K)$ is infinite dimensional and $H_2(\widehat{\lib/K})$ is of uncountable dimension.\end{theorem}

 \begin{proof} With the notation in past sections, write as above a generic element
of $\lib/K$ (resp. $\widehat{\lib/K}$)  as a sum (resp. series) $\lambda a+\sum_r\lambda_r\,b_r$ with $\lambda,\lambda_r\in\bq$. Then,
 a basis of the cycles in $\quic_2(\lib/K)$ is given by the $b_r\land b_s$ with $r<s$.  Denote
 $$
 V_p= \{ b_r\land b_s,\,\,\,\, 0\leq r<s,\, r+s= p\}$$
and note that $\dim V_{2n}= n$ and $\dim V_{2n+1}= n+1$.
Moreover, the differential on $\quic_3(\lib/K)$ satisfies
 $$d (b_r\land b_s\land a) =b_{r+1}\land b_s + b_r\land b_{s+1}$$
 so that it restricts to morphisms
 $$d\colon  V_p\land a \to V_{p+1}.$$
 Clearly each element of $\quic_2(\lib/K)$ can be written in the form $b\land b_r + d(\gamma)$ for some $b_r$ and $\gamma$. On the other hand, if $p= 2n$ is even, then
 $$b\land b_p = d\, \left( \sum_{i=0}^{n-1} (-1)^i \, b_i \land b_{2n-i-1}\, \right).$$
 It follows that $d : V_{2n-1}\land a \to V_{2n}$ is an isomorphism and that
  a basis for $H_2(\lib/K)$ is given by the elements $b\land b_r$, with $r$ odd. Therefore, $H_2(\widehat{\lib/K}) = b \land S$ where $S$ is the set of series $\sum_i \lambda_i \, b_{2i+1}$, with $\lambda_i\in \mathbb Q$. This is clearly a vector space of uncountable dimension.
\end{proof}

 All the results in this note suggest the following:

 \vspace{3mm}\noindent  {\bf Question.} Let $L$ be a finitely generated  Lie algebra such that  $\cap_n L^n= 0$, being $L=L^1\supset L^2\supset \dots$ its  central series.  If $\widehat L$ is infinite dimensional, is $H_2(\widehat L)$ always of uncountable dimension?

 \end{document}